\documentclass{article}
\usepackage{latexsym, amsxtra, amsmath, amssymb, amsthm, amscd}
\usepackage[all]{xy}

\usepackage{amssymb}
\usepackage{indentfirst}
\usepackage{amscd}
\newtheorem{theorem}{Theorem}[section]
\newtheorem{lem}[theorem]{Lemma}
\newtheorem{cor}[theorem]{Corollary}

\def\Soc{\text{Soc}}

\def\m{\frak m}

\def\q{\frak q}

\def\Mj{\frac{M}{\q_j M}}
\def\Mx{\frac{M}{xM}}
\def\Mo{\frac{M}{H^0_\m(M)}}

\def\a{\frak a}

\def\over{\overline}
\def\Hom{\text{Hom}}

\def\:{:_M}
\begin{document}\title{Asymptotic Behaviour of Parameter Ideals in Generalized Cohen-Macaulay Modules   }
\author{Nguyen Tu Cuong\footnote{Email: ntcuong@math.ac.vn} \  and Hoang Le Truong\footnote{Email: hltruong@math.ac.vn}\\
Institute of Mathematics\\
18 Hoang Quoc Viet Road, 10307 Hanoi, Viet Nam}
\date{ }
\maketitle

\begin{abstract} The purpose of this paper is  to give affirmative answers to  two  open questions as follows.  Let $(R, \m)$ be a generalized Cohen-Macaulay Noetherian local ring.  Both questions, the first question was raised by M. Rogers \cite {R} and the second one is due to S. Goto and H. Sakurai \cite {GS1}, ask whether  for every parameter ideal $\q$ contained in a high enough power of the maximal ideal $\m $ the following statements are true: (1) The index of reducibility $N_R(\q;R)$   is  independent of the choice of $\q$; and  (2)  $I^2=\q I$, where $I=\q:_R\m$.
\vspace{.2cm}\\
{\it Key words:} index of reducibility, socle,  generalized Cohen-Macaulay module, local cohomology module.\\
{\it AMS Classification:}  Primary 13H45, Secondary 13H10.
\end{abstract}

\section{Introduction}Let  $R$  be a commutative Noetherian local  ring with the maximal ideal $\m$ and residue field $\frak k=R/\m$, and let $M$ be a finitely generated $R$-module with $\dim M=d$.  Recall that a submodule of $M$ is called irreducible if it cannot be written as the intersection of two larger submodules. It is well known that every submodule $N$ of $M$ can be expressed as an irredundant intersection of irreducible submodules, and that the number of irreducible submodules appearing in such an expression depends only on $N$ and not on the expression. Thus for a parameter ideal $\q$ of $M$, the number $N_R(\q;M)$ of irreducible modules that appear in an irredundant irreducible decomposition of $\q M$ is called the index of reducibility of $\q$ on $M$. Let $N$ be an arbitrary $R$-module. We denote by Soc$(N)$ the socle of $N$. Since $\text {Soc}(N)\cong 0:_N \m \cong \Hom (\frak k, N)$
is a $\frak k$-vector space, we set $s(N) = \dim_{\frak k}\text {Soc}(N)$ the socle dimension of $N$. Then
 we have $N_R(\q;M)=s(M/\q M)$.

In 1957, D. G. Northcott \cite [Theorem 3]{N} proved that
 the index of reducibility of any parameter ideal in a Cohen-Macaulay local ring is dependent only on the ring and not on the choice of the parameter ideal. However, this property of constant index of reducibility of parameter ideals does not characterize Cohen-Macaulay modules. The first example of a non-Cohen-Macaulay Noetherian local ring having constant index of reducibility of  parameter ideals was given by S. Endo and M. Narita \cite {EN}. In 1984, S. Goto and N. Suzuki \cite {GS1} considered the supremum $r(M)$ of the index of reducibility of parameter ideals of $M$ and they showed that this number is finite provided $M$ is a generalized Cohen-Macaulay module. Recall that $M$ is said to be a generalized Cohen-Macaulay module, if  local cohomology modules $H_{\m}^i(M)$ of $M$ with respect the maximal ideal $\m$ is of finite length for $i=0,1,\ldots ,d-1$.  Moreover, they also proved that $r(M)\geqslant  \sum\limits_{i=0}^d\binom{d}{i}s(H^i_\m(R))$. Later,
S.  Goto and H. Sakurai in \cite [Corollary 3.13]{GS} showed  that
if $R$ is a Buchsbaum ring of positive dimension, then there is a
power of the maximal ideal $\m$ inside which  every parameter
ideal $\q$  has the same index of reducibility. J. C. Liu and M.
Rogers \cite{LR} refer to this by saying $R$ has eventual constant
index of reducibility of parameter ideals.  Therefore the
following question, which was raised first by M. Rogers in
\cite[Question 1.2]{R} (see also \cite [Question 1.3] {LR}), is
natural: {\it Does a generalized Cohen-Macaulay rings have
eventual constant index of reducibility of parameter ideals?}

Partial answers to this question were proved by  Rogers \cite[Theorem 1.3]{R} for a generalized Cohen-Macaulay module of dimension $d\leqslant 2$ and by Liu and  Rogers \cite[Theorem 1.4]{LR} for a generalized Cohen-Macaulay module $M$ having $H^i_{\m}(M)=0$ for all $i$ with $i\not=0,t,d$, where $t$ is some integer with  $0 <t<d$.

Our first main result in this paper is to provide a completely answer to this question.
\begin{theorem}\label{a}
Let $M$ be a generalized Cohen-Macaulay module over a Noetherian local ring $(R,\m)$ with $\dim M=d$. Then there is a positive integer $n$ such that for every parameter ideal $\q$ of $M$ contained in $\m^{n}$ the index of reducibility  $N(\q;M)$  is independent of the choice of $\q$ and is given by
$$N(\q;M)=\sum\limits_{i=0}^d\binom{d}{i}s(H^i_\m(M)).$$
\end{theorem}
In \cite{GS}, Goto and Sakurai used the study of the index of reducibility of parameter ideals in order to investigate when  the equality $I^2=\q I$ holds for a parameter ideal $\q$ of $R$, where $I=\q:\m$. Note that by results of A. Corso, C. Huneke, C. Polini and W. V. Vasconcelos \cite {CHP, CP, CPW} this equality holds for any parameter ideal in a Cohen-Macaulay local ring $R$ which is not regular or dimensional at least 2 and  $e(R)>1$, where  $e(R)$ is the multiplicity of $R$ with respect to the maximal ideal $\m$.  Goto and Sakurai generalized this and proved in \cite[Theorem 3.11]{GS} that if $R$ is a Buchsbaum ring of dimension $\dim R\geq 2$ or $\dim R= 1 $ and  $e(R)>1$,  then the equality $I^2=\q I$ holds true for all parameter ideals $\q$ contained in a high enough power of the maximal ideal $\m$. From this point of view, it is natural to ask the following question, which is due to Goto-Sakurai \cite[p. 34]{GS}:
 {\it Let   $R$ be a generalized Cohen-Macaulay ring with the multiplicity $e(R) > 1$. Is there  a positive  integer $n$ such that $I^2=\q I$ for every parameter ideal $\q$ contained in $\m^n$?}

As a consequence of Theorem \ref {a} we obtain the second main result of the paper, which is an affirmative answer  to this  question.
\begin{theorem}\label{b}
Let $R$ be a generalized Cohen-Macaulay ring and assume that $\dim R\ge 2$ or $\dim R=1$, $e(R)>1$. Then there exists a positive integer $n$ such that $I^2=\q I$ for every parameter ideal $\q\subseteq \m^n$, where $I=\q:\m$.
\end{theorem}
 Our goal for proving Theorem \ref {a}  is to show by induction on
$d=\dim M$ that there is an enough large integer $n$ such that
$N(\q;M)=\sum\limits_{i=0}^d\binom{d}{i}s(H^i_\m(M))$ for every
parameter ideal $\q\subseteq \m^n$. Therefore we give in the
Section 2 several lemmata on the asymptotic behaviour of parameter
ideals in a generalized Cohen-Macaulay module $M$ in order to
prove the following key result in Section 3 (see Theorem \ref {2a}): Let $M$ be a
generalized Cohen-Macaulay $R$-module. Then there exists a enough
large integer $k$ such that
$$s(H^i_\m(\frac{M}{(x_1,\ldots,x_{j+1})M})
=s(H^i_\m(\frac{M}{(x_1,\ldots,x_{j})M}))+s(H^{i+1}_\m(\frac{M}{(x_1,\ldots,x_{j})M})),$$ for every parameter ideal
$\q=(x_1,\ldots,x_d)\subseteq \m^k$ and for all  $0\leqslant i+j
\leqslant d-1$. The last Section is devoted to prove  the main results and their consequences.

\section{Some  auxiliary lemmata}

Throughout this paper we fix the following standard notations: Let $R$ be a Noetherian local commutative ring
with maximal ideal $\m$,  $\frak k=R/\m$ the residue field and $M$  a finitely generated $R$-module
with $\dim M=d$. Let $\q=(x_1,\ldots,x_d)$ be a parameter ideal of
module $M$.  We denote by $\q_i$  the ideal $(x_1,\ldots,x_i)R$ for $i=1,\ldots , d$  and stipulate that  $\q_0$ is  the zero ideal of $R$.

An  $R$-module $M$ is said to be a {\it generalized Cohen-Macaulay module} if 
$H^i_\m(M)$ are of finite length for all $i=0,1,\ldots , d-1$ (see \cite {CST}). This condition is equivalent to saying that there exists a  parameter ideal $\q=(x_1,\ldots,x_d)$ of $M$ such that $\q H^i_\m(\Mj)=0$ for all $0\le i+j < d$ (see \cite {T}), and such a parameter ideal was called a {\it standard parameter ideal} of $M$. It is well-known that if $M$ is a generalized Cohen-Macaulay module, then every parameter ideal of $M$ in a high enough power of the maximal ideal $\m$ is standard. The following  lemma can be easily derived from the basic properties of generalized Cohen-Macaulay modules.
\begin{lem}\label{312}
Let $M$ be a generalized Cohen-Macaulay $R$-module with $\dim M=d\ge 1$. Then there exists a positive integer $n_1$ such that  for all parameter ideals
$\q=(x_1,\ldots,x_d)$ of $M$ contained in $\m^{n_1}$ we have $\m^{n_1}H^i_\m(\Mj)=0$
 for  all $0\le i+j\le d-1$.
\end{lem}
\begin{proof}
Since $M$ is a generalized Cohen-Macaulay $R$-module, there is an
integer $l$ such that $\m^lH^i_\m(M)=0$ for all $0\le i\le d-1$. Let
$x\in \m^l$ be a parameter element of $M$. Since $\ell (0:_Mx)<\infty$, we have isomorphisms $H^i_\m(M)\cong H^i_\m(\Mx)$ for all $i\ge 1$, and so that  the
sequences
$$\xymatrix{0\ar[r]&H^i_\m(M)\ar[r]&H^i_\m(\Mx)\ar[r]&H^{i+1}_\m(M)\ar[r]&0}$$
are exact 
for all $0\le i\le d-2$. Therefore $\m^{2l}H^i_\m(\Mx)=0$ for all $0\le i\le d-2$. Now, set $n_1= 2^{d-1}l$. We can use the fact above to prove that for all parameter ideals $\q=(x_1,\ldots,x_d)$ of $M$ contained in $\m^{n_1}$ and $0\le i+j\le d-1$, it holds  $\m^{n_1}H^i_\m(\Mj)=0$. 
\end{proof}

In order to prove the next lemma, we need a result of W. V. Vasconcelos  on the reduction number of an ideal in local rings.
 Let $J$ and $K$ be two ideals of $R$ with $J\subseteq K$.
 The ideal $J$ is called a reduction of $K$ with respect to $M$ if
$K^{r+1}M=JK^rM$ for some integer $r$, and the least of such integers is
denoted by $r_J(K,M)$.  Then  the big reduction number  $\text{bigr}(K)$  of $K$ with
respect to $M$ was defined by 
$$\text{bigr}(K)=\sup\{r_J(K,M)|\ J \text{  is a reduction of } K \text{  with respect to }
M\}.$$
It is known that  there always exists a reduction ideal for any ideal $K$ provided the residue field $\frak k$ of $R$ is infinite. Especially, if $K$ is $\m$-primary then any minimal reduction ideal of $K$  with respect to $M$ is a parameter ideal of $M$. Moreover, it was shown by Vasconcelos  \cite{Va} that $\text{bigr}(K)$ is finite for any ideal $K$.

\begin{lem}\label{22}
Let $M$ be a generalized Cohen-Macaulay $R$-module with $\dim
M=d\ge 1$. Then there exists a positive  integer $n_2$ such that for all parameter ideals $\q=(x_1,\ldots,x_d)$ of
$M$ contained in $\m^{n_2}$ and $0\le j < d$ we have $$\m^{n_2}\Mj\cap
H^0_\m(\Mj)=0.$$
\end{lem}
\begin{proof}

Note first that by  the faithfully  flat homomorphism $R\to
R[X]_{\m R[X]}$ as a basic change,  we can assume without any loss of generality that the residue field $\frak k$ of $R$ is infinite.
By Lemma \ref{312} there is an integer $n_1$ such that
$H^0_\m(\Mj)=0:_{\Mj}\m^{n_1}$ for all parameter ideals
$\q$ contained in $\m^{n_1}$ and $j<d$. Set
$K=\m^{n_1}$ and $n_2=(\text{bigr}(K)+1)n_1$.  Then for any
parameter ideal $\q=(x_1,\ldots,x_d)$ of $M$ contained in $\m^{n_2}$ and any  $0\le j<d$, there is a parameter ideal $\a=(a_{j+1},\ldots,a_d)$ of $\Mj$
contained in $K$, which is a reduction of $K$ with respect to $\Mj$,  such that 
$$\a K^{r_\a(K,\Mj)} \Mj=K^{r_\a(K,\Mj)+1}\Mj.$$ 
Since $r_\a(K,\Mj)\le r_\a(K,M)\le 
\text{bigr}(K)<\infty$,
we have
$$\m^{n_2}\Mj\cap H^0_\m(\Mj)=\a K^{\text{bigr}(K)}\Mj\cap H^0_\m(\Mj)\subseteq \a\Mj\cap H^0_\m( \Mj).$$
Therefore  it is enough   to prove that $\a\Mj\cap H^0_\m( \Mj)=0$.  In fact,
let $m\in \a \Mj\cap H^0_\m(\Mj)$.  Write
$m=a_{j+1}m_{j+1}+\ldots+a_dm_d$,  where $m_i\in \Mj$ for all
$i=j+1,\ldots,d$. Since $\Mj$ is a generalized Cohen-Macaulay
module and $\a$ a standard parameter ideal of $\Mj$ by Lemma \ref{312}, we get that 
$$m_d\in
(a_{j+1},\ldots,a_{d-1})\Mj:a_d^2=(a_{j+1},\ldots,a_{d-1})\Mj:a_d.$$
 It follows that 
$$\a \Mj\cap H^0_\m(\Mj)\subseteq (a_{j+1},\ldots,a_{d-1})
\Mj\cap H^0_\m(\Mj).$$ If $j+1<d-1$, we can continue the procedure  above again so that after $(d-j)$-times we obtain
 $$\a \Mj\cap
H^0_\m(\Mj) \subseteq a_{j+1}
\Mj\cap H^0_\m(\Mj)\subseteq a_{j+1}
\Mj\cap (0:_{\Mj}a_{j+1})= 0$$
as required. 
\end{proof}

\begin{lem}\label{23}
Let $M$ be a finitely generated $R$-module with $\dim M=d\ge 1$.
Let $k$ and $\ell$ be two positive integers. Then there
exists an integer $n_3> \ell$ such that
$$(\m^{n_3}+H^0_\m(M)):\m^{k}\subseteq\m^{\ell}M+H^0_\m(M).$$

\end{lem}
\begin{proof}
Let $\over M=\Mo$. Then there is an $\over M$-regular element $a$ contained in $\m^{k}$. By the Artin-Rees Lemma, there
exists a positive integer $m$ such  that $\m^{\ell+m}\over
M\cap a\over M=\m^\ell(\m^{m}\over M\cap a\over M)$. Set $n_3=\ell+m$. We have
$$a(\m^{n_3}\over M:\m^k)\subseteq a(\m^{n_3}\over M:a)=\m^{n_3}\over M\cap a\over M=\m^{\ell}(\m^m\over M\cap a\over M),$$
so that $a(\m^{n_3}\over M:\m^k)\subseteq a\m^{\ell}
\over M$. It follows from the regularity of  $a$ that
$\m^{n_3}\over M:\m^k\subseteq \m^{\ell} \over M$.
Hence $(\m^{n_3}M+H^0_\m(M)):\m^k\subseteq \m^{\ell}
M+H^0_\m(M)$ as required.
\end{proof}

\begin{lem}\label{251a}Let $M$ be a finitely generated $R$-module with $\dim M=d\ge 1$.
Then there exists a positive integer $n_4$ such that  for all ideals $K\subseteq \m^{n_4}$ we have
 $$(KM+H^0_\m(M)):\m =KM:\m+H^0_\m(M).$$
\end{lem}
\begin{proof}
Since $H^0_\m(M)$ have finite length, there exists an integer $\ell$ such that
$\m^\ell M\cap H^0_\m(M)=0$. By Lemma \ref{23}, there is an integer $n_4> \ell$ such
that for all ideals $K\subseteq\m^{n_4}$  we have
$$(KM+H^0_\m(M)):\m\subseteq (\m^{n_4}M+H^0_\m(M)):\m\subseteq \m^{\ell}M+H^0_\m(M).$$ 
 Let $b\in(KM+H^0_\m(M)):\m$. Write $b=\alpha+\beta$ with $\alpha\in \m^{\ell}M$ and $\beta\in H^0_\m(M)$.
 Then,  since $K\subseteq \m^{n_4}\subseteq\m^{\ell+1}$,  $$\m\alpha\subseteq \m^{\ell+1}M\cap(KM+H^0_\m(M))=KM+\m^{\ell+1}M\cap H^0_\m(M)=KM.$$
  Thus $\alpha\in KM:\m$ and so that $(KM+H^0_\m(M)):\m
 =KM:\m+H^0_\m(M)$. 
\end{proof}

\begin{lem}\label{251}Let $M$ be a generalized Cohen-Macaulay $R$-module with $\dim M=d\ge 1$.
Then there exists a positive integer $n_5$ such that for all parameter ideals $\q=(x_1,\ldots,x_d)$ of
$M$ contained in $\m^{n_5}$ and $0\le j<i\le d$ we have
$$[\frac{\q_iM}{\q_jM}+H^0_\m( \frac{M}{\q_jM})]:\m
=\frac{\q_iM}{\q_jM}:\m+H^0_\m(\frac{M}{\q_jM}).$$ 
.
\end{lem}
\begin{proof}
Let $n_1$ and $n_2$ be two integers as in Lemma \ref{312} and Lemma \ref{22}, respectively. 
By Lemma \ref{23}, there  always exists an integer
$n_5>n_2$ such that $ (\m^{n_5}M+H^0_\m(M)):\m^{n_1 +1}\subseteq
\m^{n_2} M+H^0_\m(M)$ . Let $\q=(x_1,\ldots,x_d)$ be a parameter ideal of $M$ contained in $\m^{n_5}$. For all $0\le j<i\le d$, we have $H^0_\m(\Mj)=0:_{\Mj}\m^{n_1}$ by Lemma \ref{312},  and so that
$$\begin{aligned}
(\frac{\q_i M}{\q_j M}+H^0_\m( \Mj)):\m&\subseteq \frac{\m^{n_5} M}{\q_j M}:\m^{n_1+1}\\
&= \frac{\m^{n_5} M:\m^{n_1+1}}{\q_j M}\subseteq
\frac{\m^{n_2} M}{\q_j M}+H^0_\m(\Mj).\end{aligned}$$
 Let $b\in(\frac{\q_iM}{\q_jM}+H^0_\m( \Mj)) :\m$. Write $b=\alpha+\beta$ with $\alpha\in \frac{\m^{n_2} M}{\q_jM}$ and $\beta\in
H^0_\m(\Mj)$. Since $\q_i\subseteq \m^{n_5}\subseteq \m^{n_2+1}$, we get by Lemma
\ref{22} that $$\m\alpha\subseteq \frac{\m^{n_2 +1}M}{\q_jM}\cap(\frac{\q_iM}{\q_jM}+H^0_\m(\Mj))=\frac{\q_iM}{\q_jM}+ \frac{\m^{n_2+1}M}{\q_jM}\cap H^0_\m( \Mj)=\frac{\q_iM}{\q_jM}.$$ Therefore $\alpha\in \frac{\q_iM}{\q_jM}:\m$, and so that
$$(\frac{\q_iM}{\q_jM}+H^0_\m( \Mj)):\m =\frac{\q_iM}{\q_jM}:\m+H^0_\m(\Mj)$$
as required.
\end{proof}

\section{The socle dimension of local cohomology modules}
Let $\q=(x_1,\ldots,x_d)$ be a parameter ideal of the module $M$. For
each positive integer $n$, we denote by $\q(n)$ the ideal
$(x_1^n,\ldots,x_d^n)$. Let $K_*(\q(n))$ be the Koszul complex of
$R$ with respect to the ideal $\q(n)$ and
$$H^*(\q(n);M)=H^*(\Hom(K_*(\q(n),M))$$ the Koszul cohomology
module of $M$. Then  the family
$\{H^i(\q(n);M)\}_{n\ge1}$ naturally forms an inductive system of
$R$-modules for every $i\in\Bbb Z$, whose inductive limit is just the $i$-th local cohomology module
$$H^i_\m(M)=H^i_\q(M)=\varinjlim\limits_nH^i(\q(n);M).$$ 
The following result is due to Goto and Suzuki.
\begin{lem}[\cite{GS1}, Lemma 1.7]\label{split}  Let  $M$  be a finitely generated $R$-module, $x$ an $M$-regular element and $\q=(x_1,\ldots,x_r)$  an ideal of $R$ with $x_1= x$. Then there exists a splitting  exact sequence for each $i \in\Bbb Z$,
$$0\to H^i(\q;M)\to H^i(\q;\Mx)\to H^{i+1}(\q;M)\to0.$$
\end{lem}
The next result is due to Goto and Sakurai.
\begin{lem}[\cite{GS} Lemma 3.12]\label{sur}
Let $R$ be a Noetherian local ring with the maximal ideal $\m$ and
$r = \dim R \ge1$. Let $M$ be a finitely generated $R$-module.
Then there exists a positive integer $\ell $ such that for all parameter ideals $\q=(x_1,\ldots,x_d)$  of $M$ contained in $\m^\ell$ and
all  $i\in\Bbb Z$, the canonical homomorphisms on socles
$$\rm{ Soc }(H^i(\q,M))\to \rm{ Soc }(H^i_\m(M))$$
are surjective.
\end{lem}

The following theorem is the key to proofs of  main results of the paper.
\begin{theorem}\label{2a}
Let $M$ be a generalized Cohen-Macaulay $R$-module with $\dim
M=d\ge 1$. There  there exists a positive integer $k$ such that for all parameter ideal $\q$ of
$M$ contained in $\m^k$ and  $d> i+j \ge 0$ we have
$$s(H^i_\m(\frac{M}{\q_{j+1}M})) =s(H^i_\m(\Mj))+s(H^{i+1}_\m(\Mj)),$$
where $s(N) = \dim_{\frak k}\rm {Soc }(N)$ the socle dimension of the $R$-module $N$.
\end{theorem}
\begin{proof}
We set $k=\max\{n_1, n_2, n_5 , \ell\}+1,$ where $n_1,$ $n_2$, $n_5 $ and $\ell$ are integers as in Lemma \ref{312}, \ref{22} , \ref{251}, and \ref{sur}, respectively. It  will be shown that this integer $k$ is just the required integer of the theorem. Let $\q=(x_1,\ldots,x_d)$ be a parameter ideal of $M$ contained in $\m^k$. We denote by $M_j$  the module $\Mj$ and $\over M_j$  the module $\frac{M_j}{H^0_\m(M_j)}$. It should be noted here that $M_j$ and $\over M_j$ are generalized Cohen-Macaulay modules having $(x_{j+1},\ldots , x_d)$ as a standard parameter ideal by Lemma \ref{312}. Then the proof of Theorem \ref{2a} is divided into two cases.
\vskip 0.3cm
\noindent
{\it First case}: $i=0$.
Because of the choose of $k$, the ideal $\q$ is a standard parameter ideal of $M$ and so that $x_{j+1}H^1_\m(\over M_j)=0$ for all $0\le j<d$. Thus we have $$ H^1_\m(M_j)\cong H^1_\m(\over M_j)\cong H^0_\m(\frac{\over M_j}{x_{j+1}\over M_j}).$$ Therefore, we
get  by Lemma \ref{251} that 
$$\begin{aligned}s(H^1_\m(M_j))&=s(H^0_\m(\frac{\over M_j}{x_{j+1}\over M_j}))
=\ell(\frac{(\q_{j+1}M_j+H^0_\m(M_j)):\m}{\q_{j+1}M_j+H^0_\m(M_j)})
\\&=\ell(\frac{\q_{j+1}M_j:\m+H^0_\m(M_j)}{\q_{j+1}M_j+H^0_\m(M_j)})
\\&=\ell(\frac{\q_{j+1}M_j:\m}{(\q_{j+1}M_j:\m)\cap(\q_{j+1}M_j+H^0_\m(M_j))})
\\&=\ell(\frac{\q_{j+1}M_j:\m}{\q_{j+1}M_j+(\q_{j+1}M_j:\m)\cap H^0_\m(M_j)}).
\end{aligned}$$
 Let $a\in(\q_{j+1}M_j:\m)\cap H^0_\m(M_j)$. We
see by Lemma \ref{22} that 
$$\m a\in\q_{j+1}M_j\cap H^0_\m(M_j)=0.$$ Therefore
$(\q_{j+1}M_j:\m)\cap H^0_\m(M_j)=0:_{M_j}\m$, and so that
$$\begin{aligned}s(H^1_\m(M_j))&=\ell(\frac{\q_{j+1}M_j:\m}{\q_{j+1}M_j+0:_{M_j}\m})\\
&=\ell(\frac{\q_{j+1}M_j:\m}{\q_{j+1}M_j})-\ell(\frac{\q_{j+1}M_j+0:_{M_j}\m}{\q_{j+1}M_j})\\
&=\ell(\frac{\q_{j+1}M_j:\m}{\q_{j+1}M_j})-\ell(0:_{M_j}\m)\\
&=s(H^0_\m(M_{j+1}))-s(H^0_\m(M_j)).
\end{aligned}$$
Hence, we have
$s(H^0_\m(M_{j+1}))= s(H^0_\m(M_j))+s(H^1_\m(M_j))$ for all $0\le j<d$.
\vskip 0.3cm
\noindent
{\it Second case}: $i\ge 1$. We first claim by induction on $j$ that for all $i\ge 1$ and $d> i+j \ge 1$, the canonical homomorphisms on socles
$$\alpha _j^i:\Soc(H^i(\q,\over M_j))\to \Soc(H^i_\m(\over M_j))$$
are surjective. For the case $j=0$, we consider the following commutative diagram
$$\xymatrix{H^i(\q;M)\ar[d]^{f _i}\ar[r]&H^i(\q,\over M_0)\ar[d]^{g_i}\\
H^i_\m(M)\ar[r]^{\pi _i}&H^i_\m(\over M_0),}$$ 
where $\pi_i$ are isomorphisms  for all $i\ge 1$. By Lemma \ref{sur},  the   homomorphism $f_i$ induces a surjective homomorphism  $\Soc(H^i(\q,M))\to \Soc(H^i_\m( M))$ on the socles.  Therefore we get by applying the functor $\Hom(\frak k,*)$ to the diagram above that 
$$\alpha _0^i:\Soc(H^i(\q,\over M_0))\to \Soc(H^i_\m(\over M_0))$$
are surjective for all $i\ge 1$.
 Now assume that $j\ge 1$.  Since $(x_{j+1},\ldots , x_d)$ is a standard parameter ideal of $\over M_j$ and  $x_{j+1}$  an $\over M_j$-regular element, we have for all $d> i+j \ge 1$ the following commutative diagram
$$\xymatrix{0\ar[r]&H^i(\q;\over M_j)\ar[d]\ar[r]&H^i(\q;\frac{\over M_j}{x_{j+1}\over M_j})\ar[r]\ar[d]&H^{i+1}(\q;\over M_j)\ar[r]\ar[d]&0\\
0\ar[r]&H^i_\m(\over M_j)\ar[r]&H^i_\m(\frac{\over M_j}{x_{j+1}\over M_j})\ar[r]&H^{i+1}_\m(\over M_j)\ar[r]&0}$$
with exact rows, where the upper row is split exact by Lemma \ref{split}.  Therefore, by applying the functor $\Hom(\frak k,*)$, we obtain for all $d> i+j \ge 1$
the commutative diagram
$$\xymatrix{0\to \Soc(H^i(\q;\over M_j))\ar[d]^{\alpha _j^i}\ar[r]&\Soc(H^i(\q;\frac{\over M_j}{x_{j+1}\over M_j}))\ar[r]\ar[d]^{\beta_{j+1}^i}&\Soc(H^{i+1}(\q;\over M_j))\to0\ar[d]^{\alpha _j^{i+1}}\\
0\to
\Soc(H^i_\m(\over M_j))\ar[r]&\Soc(H^i_\m(\frac{\over M_j}{x_{j+1}\over M_j}))\ar[r]&\Soc(H^{i+1}_\m(\over M_j))}$$
with exact rows.  By the inductive hypothesis,  the homomorphisms $\alpha _j^i$ and $\alpha _j^{i+1}$ are surjective for all $i\ge1$. Thus
$\beta _{j+1}^i$ are surjective for all $i\ge 1$. Since $\over M_j$ is generalized Cohen-Macaulay, it is easy to check that $H^i_\m(\frac{\over M_j}{x_{j+1}\over M_j})\cong H^i_\m(\over M_{j+1})$ for all $i\ge1$. It follows from the commutative diagram
$$\xymatrix{\Soc H^i(\q;\frac{\over M_j}{x_{j+1}\over M_j})\ar[d]^{\beta_{j+1}^i}\ar[r]&\Soc H^i(\q,\over M_{j+1})\ar[d]^{\alpha_{j+1}^i}\\
\Soc H^i_\m(\frac{\over M_j}{x_{j+1}\over M_j})\ar[r]^{\cong}&\Soc H^i_\m(\over M_{j+1})}$$ 
that
the homomorphism $\alpha _{j+1}^i:\Soc(H^i(\q,\over M_{j+1}))\to \Soc(H^i_\m(\over M_{j+1}))$ are surjective for all $d>i+j\ge 1$, and the claim is proved. Next, from the proof of the claim 
we obtain exact sequences
$$\xymatrix{0\to
\Soc(H^i_\m(\over M_j))\ar[r]&\Soc(H^i_\m(\frac{\over M_j}{x_{j+1}\over M_j}))\ar[r]&\Soc(H^{i+1}_\m(\over M_j))\to0},$$
and so that $s(H^i_\m(\frac{\over M_j}{x_{j+1}\over M_j})) =s(H^i_\m(\over M_j))+s(H^{i+1}_\m(\over M_j))$ for all $i\ge 1$ and $d>i+j\ge 0$. Therefore, since 
 $H^i_\m(\over M_j)\cong H^i_\m(\Mj)$ for all $i\ge 1$, we have $$s(H^i_\m(\frac{M}{\q_{j+1}M})) =s(H^i_\m(\Mj))+s(H^{i+1}_\m(\Mj))$$
 for all $i\ge 1$ and $d>i+j\ge 1$, and the proof of Theorem \ref {2a} is complete.
\end{proof}

\section{Proofs of main results}
Theorem \ref{a} is now an easy consequence of Theorem \ref{2a}.
\begin{proof}[Proof of Theorem \ref{a}]
By virtue of Theorem \ref{2a} we can show by induction on $d$ that there exists an integer $n$ such that for
every parameter ideal $\q=(x_1,\ldots,x_d)$ of $M$ contained in
$\m^{n}$ 
we have
$$N(\q;M)=s(H^0_\m(\frac{M}{\q M}))=\sum\limits_{i=0}^d\binom{d}{i}s(H^i_\m(M)).$$
\end{proof}

\begin{cor}\label{sup}
Let $M$ be a generalized Cohen-Macaulay $R$-module. Then
$$\sup\{N(\q;M)| \q\text{ is a standard parameter ideal of }M\}=\sum\limits_{i=0}^d\binom{d}{i}s(H^i_\m(M)).$$
\end{cor}
\begin{proof} Let  $\q=(x_1,\ldots,x_d)$ be a standard parameter ideal of $M$. By basic properties of the theory of generalized Cohen-Macaulay modules we can show by induction on $t$ that 
$$s(H^i_\m(\frac{M}{(x_1,\ldots,x_t)M}))\le \sum\limits_{j=0}^{t}\binom{t}{j}s(H^{j+i}_\m(M)).$$
for all $d \ge i+t\ge 0$. Therefore  the Corollary follows by the inequality above in the case $t=d$, $i=0$ and Theorem \ref{a}.
\end{proof}

In the rest of this paper, we denote
$$S(M)=\sum\limits_{i=0}^d\binom{d}{i}s(H^i_\m(M)).$$

 \begin{proof} [Proof of Theorem \ref{b}]
Let $n=\max \{n_1, n_4, k\}$, where $n_4$ and $k$ are integers in Lemma \ref{312}, Lemma \ref{251a}  and Theorem \ref{2a} (for the case $M=R$), respectively. We will prove that $I^2=\q I$ for all parameter ideals $\q=(x_1,\ldots,x_d)$ of $R$ contained in $\m^n$, where $I=\q:_R\m$. Let $\dim R=d$ and $ \over R=\frac{R}{H^0_\m(R)}$.  Then by Lemma
\ref{251a} we have
$$(\q+H^0_\m(R)):_R\m =\q:_R\m+H^0_\m(R),$$ and so that  $I\over R=\q \over R:\m
\over R$. 
\vskip0.3cm
\noindent
{\it Case 1:  $e(R)=1$ and $d\ge 2$}. Since $\over R$ is unmixed, it is well-known in this case hat
 $\over R$ is a regular local ring of dimension $d\ge 2$. We have
  $(I\over R)^2=\q \over R I\over
R$ by
Theorem 2.1 in \cite{CPW}. Therefore $I^2\subseteq \q I+H^0_\m(R)$ and so that
$I^2\subseteq\q I+I^2\cap H^0_\m(R)$. But, $I^2\cap H^0_\m(R)\subseteq \q\cap H^0_\m(R)=0$ by
Lemma \ref{22}. Thus  $I^2=\q I$ in this case.
\vskip0.3cm
\noindent
{\it Case 2:  $e(R)>1$}. By the choose of $n$, the parameter ideal $\q$ is standard Lemma \ref{312}  
and $N(\q;R)=S(R)$ by Theorem
\ref{a}. Thus, it is enough for us to prove  that if $N(\q;R)=S(R)$ for some standard
parameter ideal $\q=(x_1,\ldots,x_d)$ of $R$ contained in $\m^n$ 
then $I^2=\q I$.  Indeed, we argue by induction
 $d$. Let  $d=1$. Then $\over R$ is a non-regular Cohen-Macaulay ring, and the conclusion follows with the same method as used in the proof of case 1.  Now assume that $d\ge 2$.  Set $R^\prime=\frac{R}{(x_1)}$.
By Theorem \ref{2a}, we have $S(R)=S(R^\prime)$, and so that $N(\q R^\prime;
R^\prime)=S(R^\prime)$.  Therefore $ (IR^\prime)^2=\q R^\prime IR^\prime$ by the inductive
hypothesis. It follows that  
  $I^2\subseteq
(x_2,\ldots,x_d)I+(x_1)$, and so that
$I^2\subseteq(x_2,\ldots,x_d)I+(x_1)\cap I^2$. Let $a\in(x_1)\cap
I^2$ and we write $a=x_1b$ with $b\in R$. Since $e(R)>1$, by
Proposition (2.3) in \cite{GS}, we have $\m I^2=\m \q^2$.
Therefore $\m a=x_1\m b\subseteq (x_1)\cap \q^2$. Since the
parameter ideal $\q$  is standard, 
$(x_1)\cap\q^2=x_1\q$ and $H^0_\m(M)=0:_Rx_1$. Thus $\m b\subseteq
(x_1\q):_Rx_1=\q+0:_Rx_1$, and so that $b\in
(\q+0:_Rx_1):_R\m=\q:_R\m+0:_Rx_1$ by Lemma \ref{251a}. Therefore $a\in x_1I$, and so that $(x_1)\cap
I^2=x_1I$. Hence $I^2=\q I$ as required.
\end{proof}

\begin{cor}
Let $R$ be a generalized Cohen-Macaulay local ring with
multiplicity $e(R)>1$. Then for sufficiently large $n$, we have
$$\mu(I)= d+S(R)$$
 for all parameter ideals $\q$ contain in $\m^n$,
where $\mu(I)$ is the minimal number of generators of the ideal
$I=\q:\m$.
\end{cor}
\begin{proof}
Choose the integer $n$ as in Theorem \ref{a} (for the case $M=R$). Then  
$$\frac{I}{\q}\cong \Hom(\frak{k},\frac{R}{\q})\cong\frak{k}^{S(R)}$$ 
by Theorem
\ref{a}. Since $e(R)>1$, by Proposition 2.3 in
\cite{GS}, we get that $\m I=\m\q$. Therefore
$$\mu(I)=\ell(\frac{I}{\m I})=\ell(\frac{I}{\m\q})=\ell(\frac{I}{\q})+\ell(\frac{\q}{\m\q})=S(R)+d$$
 as required.
\end{proof}

\end{document}